\documentclass[12pt]{article}
\usepackage{amssymb,latexsym}
\usepackage{amsmath}
\usepackage{tikz}
\usepackage{graphicx}
\usepackage[all,dvips]{xy}
\newtheorem{teorema}{Theorem}[section]

\newtheorem{lema}[teorema]{Lemma}
\newtheorem{proposicao}[teorema]{Proposition}

\newenvironment{prova}{\setlength{\parindent}{0pt}\textbf{Proof.}}{\hspace{\stretch{1}}$\Box$}

\begin{document}

\title{\textbf{Paradeduction in axiomatic formal systems}}

\author{Edelcio G. de Souza$^{1}$ \\ 
Alexandre Costa-Leite$^{2}$ \\
Diogo H.B. Dias$^{3}$ \\
\footnotesize $^{1, 3}$ University of Sao Paulo (BR) \ $^{2}$ University of Brasilia (BR)  \\
\footnotesize $^{1}$ edelcio.souza@usp.br \ $^{2}$costaleite@unb.br \ $^{3}$diogo.bispo.dias@gmail.com}

\date{}

\maketitle

\begin{abstract}
This paper presents the concept of \emph{paradeduction} in order to justify that we can overlook
contradictory information taking into account only what is consistent. Besides
that, it uses paradeduction to show that there is a way to transform any logic, introduced
as an axiomatic formal system, into a paraconsistent one.

\end{abstract} 
\vspace{1.5cm}
\section{Introduction}

The general motivation of this paper lies in the role consistent sets play in science, philosophy and reasoning in general. In ordinary life, we are recipients of much contradictory information, and sometimes our convictions and theories are also contradictory. How do we deal with them? The main \emph{motto} of this paper is that we can ignore pieces of contradictory information and concentrate only on consistent parts of it: always reason with consistent sets of information. This is indeed what occurs in Courts of Law when a judge produces a decision or in philosophy when someone reasons with contradictory viewpoints. In this sense,  we try to provide a formalization of a practice already common when we need to deal with inconsistent data. For instance, from Bohr's atomic model we  infer that electrons can orbit an atom's nucleus without any emission of energy. But if we take into account a regular component of Bohr's model such as Maxwell's equations, then we know that an orbiting electron emits energy (for more on this topic, cf. \cite{brown2}). In this scenario, the calculations were made by informally splitting the information into consistent subsets. The important thing to note is that in many different situations, when facing an inconsistent set of information, we intuitively deal with consistent subsets of the original one.  Here, we intend to make this procedure explicit by formalizing it.

Some comments on previous work dealing with consistent sets seems to be useful as there is a vast literature on using consistent subsets for reasoning with inconsistent information. The general procedure is the following: given an inconsistent set of premises, criteria to single out one or more of its consistent subsets are established. Articles \cite{avron}, \cite{benferhat}, \cite{vandeputte} and \cite{subrahmanian} use some ranking system between consistent subsets to determine preferred bases for inferences, \cite{amgoud} uses maximal consistent subsets of knowledge bases while \cite{brown} defines chunks of consistent information from a given inconsistent set. Also, a specific logic is applied to the selected sets: \cite{avron}, \cite{brown} and \cite{subrahmanian} use classical logic, \cite{amgoud} uses an abstract monotonic logic, \cite{avron} and \cite{vandeputte} explore reasoning with inconsistent set using non-monotonic logics. There are a few cases of generalization for different logics: \cite{arieli2} studies sequent-based inferences in non-classical Tarskian logics, and \cite{brown} mentions the possibility of generalizing their work to non-classical logics. 

Despite the fact that this paper is, in some sense, included in this tradition, our approach seems to go in a different direction, given that our starting point is not an inconsistent set, but a certain logic in which inconsistency can arise. In this scenario our main focus is not on reasoning with consistent subsets of an inconsistent set of information, but on the paraconsistentization of the original logic, i.e., on the process of turning an explosive logic into a paraconsistent one, while keeping some properties from the initial logic. In this sense, our work is not only distinct from the available literature, but complementary to some of them. For instance, \cite{benferhat} claims that the mechanism for handling inconsistent knowledge bases should be strictly syntactic, since they are trivially semantically equivalent, whereas our method also makes sense of a semantically inconsistent set. Additionally, \cite{arieli} proposes to establish criteria for an ideal logic for reasoning with inconsistency by imposing the preservation of as much as possible from classical logic as a constraint, this is a strong demand for our purposes.

As usual, we take a \emph{consequence structure} as a pair $(X,Cn)$ such that $X$ is a non-empty set and $Cn$ is an operation in the powerset of $X$: 
$$
Cn: \wp(X) \rightarrow \wp(X)
$$
In this way, if $A \subseteq X$, then $Cn(A)$ is the set of \emph{consequences} of $A$. Moreover, if $Cn(A) \neq X$, we say that $A$ is $Cn$-\emph{consistent}. Exploring logics from this abstract perspective is the basic idea developed by A. Tarski in \cite{tarski1} and \cite{tarski2}. 

Dealing with consequence structures (i.e. logics) in a recent paper (see \cite{desouza}), a method is proposed to show that it is possible to turn any explosive logic into a system able to deal with contradictions, answering, therefore, a question proposed in \cite{costaleite}: how to convert a given explosive logic into a paraconsistent one? In \cite{desouza}, concepts from category theory have been used to define an endofunctor in the category of consequence structures in such a way that given a logic in which \emph{ex falso} holds, it is possible to generate a paraconsistent counterpart of it by means of a \emph{paraconsistentization functor}. Indeed, it is showed that, given a consequence structure $(X,Cn)$, it is feasible to construct a paraconsistent version of it, denoted by $(X,Cn_{P})$, in such a way that, for $A \subseteq X$ and $a \in X$, we have that: $a \in Cn_{P}(A)$ if, and only if, there is a $A' \subseteq A$ $Cn$-consistent such that $a \in Cn(A')$. This is highly abstract.  

In this paper, instead of examining categories of logics, we consider something less general, i.e. axiomatic formal systems, in order to deal with the following problem: given an explosive (axiomatic) formal system, how to turn it into a paraconsistent one? Specifically, we start using formal systems, in Hilbertian style, with their standard notion of deduction. Then, we introduce the concept of \emph{paradeduction} (section 4) and show that it is related with the construction proposed in \cite{desouza}: paradeducibility in axiomatic formal systems corresponds in some sense to a paraconsistentization functor in the category of consequence structures (Proposition 4.2). Then, we are able to determine whether some properties are invariant under paradeduction (theorem 4.3) providing, thus, an answer to a question left open in \cite{desouza}.
Finally, we state a sufficient condition for the notion of \emph{paradeduction} to be paraconsistent (theorem 5.1).

\section{Axiomatic formal systems}

From the viewpoint of proof theory, logics can be introduced in some different ways. A logician
could use Gentzen-style presentations such as natural deduction or sequent calculi. Alternatively, logics could be proof-theoretically introduced by methods such as resolution or \emph{tableaux},
or even any other method still waiting to be developed could be used. In this paper, we take into consideration
logics developed by means of axiomatic systems (i.e. a Hilbert-style systems). These are
\emph{axiomatic formal systems}. In order to proceed from this syntactical perspective, we
begin with somewhat standard terminology in the realm of proof theory.

Let $n$ be a positive integer and $X$ a (non-empty) set. We denote by 
$\wp^{(n)}(X)$, the set of subsets of $X$ with cardinal $n$, i.e.,
$$
\wp^{(n)}(X) := \{A \subseteq X: |A| = n \}.
$$

An \emph{inference rule of degree $n$ on $X$} is a binary relation $R$ such that
$$
R \subseteq \wp^{(n)}(X) \times X.
$$

If the pair $(\{x_{1},...,x_{n}\},x) \in R$, we say that $x$ is an
\emph{immediate consequence of $\{x_{1},...,x_{n}\}$  in virtue of the application of the inference rule $R$}. In this case, we can use the following notation: 
$$
\frac{x_{1},...,x_{n}}{x}R.
$$

Let $\mathbb{R}$ be a family of inference rules on $X$. We say that $x$ is an \emph{immediate
consequence} of $\{x_{1},...,x_{n}\}$  if $x$ is an immediate consequence of  $\{x_{1},...,x_{n}\}$ in virtue of the application of an inference rule $R \in \mathbb{R}$.

An \emph{axiomatic formal system} $S$ is a triple $S = (X,\mathbb{A},\mathbb{R})$ in which $X$ is a non-empty set whose elements are called \emph{formulas} of $S$; $\mathbb{A}$ is a subset of $X$ whose
elements are called \emph{axioms} of $S$; and $\mathbb{R}$ is a finite family of inference rules on $X$.

Consider a formal system $S = (X,\mathbb{A},\mathbb{R})$ and let $A \subseteq X$. A finite sequence $a_{1},...,a_{n}$ of elements of $X$ is called a $S$-\emph{deduction} (or a \emph{deduction} in $S$) from $A$,
if for each $1 \leq i \leq n$, it is the case that:

\begin{enumerate}
	\item $a_{i} \in A$ or;
	\item $a_{i} \in \mathbb{A}$ or;
	\item $a_{i}$ is a consequence of precedent formulas in the sequence by the application of an inference rule.
\end{enumerate}

We say that a formula $a \in X$ is $S$-\emph{deducible} from a set $A$ of formulas if there is a $S$-deduction, $a_{1},...,a_{n}$, from $A$, such that $a = a_{n}$. To indicate that $a$ is $S$-deductible from $A$ we use the following usual notation:
$$
A \vdash_{S} a.
$$

Moreover, we say that a sequence $a_{1},...,a_{n}$ is a $S$-\emph{deduction} of $a$ from $A$.

If $A$ is a subset of $X$, we denote by $Cn_{S}(A)$ the set of all the elements of $X$ which are $S$-deducible of $A$, i.e:
$$
Cn_{S}(A) := \{a \in X:A \vdash_{S} a \}.
$$

We say that $A \subseteq X$ is $S$-\emph{consistent} if and only if $Cn_{S}(A) \neq X$. Otherwise, $A$ is said $S$-\emph{inconsistent}. We denote by $CON_{S}$ the set of all $S$-consistent subsets of $X$. In our axiomatic formal systems, we always suppose that $\emptyset$ is $S$-consistent (or, equivalently, there is at least one $S$-consistent set). We say also that $A$ is a $S$-\emph{theory} if and only if $Cn_{S}(A) = A$. $THE_{S}$ denotes the set of all $S$-theories and $THE^{*}_{S}$ refers to the set $THE_{S} \cap CON_{S}$ of $S$-consistent theories.

It is well known that the following properties hold:

\medskip
(I) If $a \in A$, then $A \vdash_{S}a$;

(II) If $A \subseteq B$ and $A \vdash_{S}a$, then $B \vdash_{S}a$;

(III)  If $B \vdash_{S}a$ and for all $b \in B$, $A \vdash_{S}b$, then $A \vdash_{S}a$;

(IV) $A \subseteq Cn_{S}(A)$;

(V) If $A \subseteq B$, then $Cn_{S}(A) \subseteq Cn_{S}(B)$;

(VI) $Cn_{S}(Cn_{S}(A)) = Cn_{S}(A)$;

(VII) If $A$ is $S$-consistent and $B \subseteq A$, then $B$ is $S$-consistent;

(VIII) If $A \vdash_{S}a$, then there exists a finite $A' \subseteq A$ such that $A'\vdash_{S}a$.

\medskip
The notion described by $CON_{S}$ is eminently syntactical considering that it deals with provability (i.e., deducibility). A formula is deducible - in a given formal system - from a set of formulas by the application of inference rules. As it is widely known, logics, in general, have two important sides, which are indeed complementary, one of them is the proof-theoretical, the other one is the semantical, which we consider in the next section.

\section{Valuation structures}

From the model-theoretical viewpoint, logics can also be introduced in very different ways. A modal logician
could use Kripke semantics to establish truth-conditions for modal operators. Or an intuitionistic logician 
could use the BHK-interpretation to analyze logical connectives. Many-valued matrices are also at our disposal to provide semantics for a given formal language. No matter which semantical technology is used, what is important is to be able to interpret a given language. Semantically, we use \emph{valuation structures}. We begin with some standard terminology (see, for instance, \cite{loparic}).

A \emph{valuation structure} is a pair $(X,\mathcal{V})$ in which $X$ is a non-empty set and $\mathcal{V}$ is a family of functions of the form:
$$
v: X \rightarrow \{0,1\},
$$
which does not contain the constant function $\mathbf{1}: X \rightarrow \{0,1\}$, such that $\mathbf{1}(x) = 1$ for all $x \in X$. Elements of $\mathcal{V}$ are called \emph{valuations} for $X$.

Let $A \subseteq X$. The set of $\mathcal{V}$-\emph{models} of $A$, denoted by $Mod_{\mathcal{V}}(A)$, is given by:
$$
Mod_{\mathcal{V}}(A) := \{v \in \mathcal{V}: v(a) = 1 \text{ for all } a \in A \}.
$$

In a similar way, if $a \in X$, then the set of $\mathcal{V}$-models of $a$ is given
by:
$$
Mod_{\mathcal{V}}(a) := Mod_{\mathcal{V}}(\{a\}) = \{v \in \mathcal{V}: v(a) = 1\}\footnote{When $v(a)=1$, we say that $a$ is \emph{true under the valuation} $v$. When that is not the case, we say that $a$ is \emph{false under the valuation $v$}.}.
$$

Note that $Mod_{\mathcal{V}}(\emptyset) = \mathcal{V}$ and $Mod_{\mathcal{V}}(X) = \emptyset$. Furthermore, if $A \neq \emptyset$, then $Mod_{\mathcal{V}}(A) = \bigcap_{a \in A}Mod_{\mathcal{V}}(a)$. It also holds that if $A \subseteq B$, then $Mod_{\mathcal{V}}(B) \subseteq Mod_{\mathcal{V}}(A)$.

A subset $A \subseteq X$ is called $\mathcal{V}$-\emph{satisfiable} if and only if $Mod_{\mathcal{V}}(A) \neq \emptyset$, i.e, there is a $v \in \mathcal{V}$ such that $v(a) = 1$, for all $a \in A$. Otherwise, $A$ is called $\mathcal{V}$-\emph{unsatisfiable}. Consider, now, a valuation structure $(X,\mathcal{V})$ and let $A \subseteq X$ and $a \in X$. We say that $a$ is a $\mathcal{V}$-\emph{consequence} of $A$ if and only if all $\mathcal{V}$-model of $A$ is also a $\mathcal{V}$-model of $a$. In this case, we use the notation $A \models_{\mathcal{V}}a$. Therefore, we have:
$$
A \models_{\mathcal{V}}a \Leftrightarrow Mod_{\mathcal{V}}(A) \subseteq Mod_{\mathcal{V}}(a).
$$

In the same way as we have done for formal systems, we define the set of $\mathcal{V}$-consequences of a subset $A \subseteq X$ as:
$$
Cn_{\mathcal{V}}(A) := \{a \in X: A \models_{\mathcal{V}}a \}.
$$

It is also well known that for valuation structures the following properties hold:

\medskip
(I) If $a \in A$, then $A \models_{\mathcal{V}}a$;

(II) If $A \subseteq B$ and $A \models_{\mathcal{V}}a$, then $B \models_{\mathcal{V}}a$;

(III) If $B \models_{\mathcal{V}}a$ and for all $b \in B$, $A \models_{\mathcal{V}}b$; then $A \models_{\mathcal{V}}a$;

(IV) $A \subseteq Cn_{\mathcal{V}}(A)$;

(V) If $A \subseteq B$, then $Cn_{\mathcal{V}}(A) \subseteq Cn_{\mathcal{V}}(B)$;

(VI) $Cn_{\mathcal{V}}(Cn_{\mathcal{V}}(A)) = Cn_{\mathcal{V}}(A)$;

(VII) If $A$ is $\mathcal{V}$-satisfiable and $B \subseteq A$, then $B$ is $\mathcal{V}$-satisfiable.

\medskip

Now, there are some known useful connections between axiomatic formal systems and valuation structures which
are important for our purposes here.

Let $S = (X,\mathbb{A},\mathbb{R})$ be a formal system and consider a valuation structure $(X,\mathcal{V})$. Note that the carrier set of a valuation structure is precisely the set $X$ of $S$ formulas. We say that the valuation structure $(X,\mathcal{V})$ is \emph{sound} with respect to the formal system $S$ if and only if for all $A \subseteq X$ and $a \in X$ we have that: if $A \vdash_{S}a$, then $A \models_{\mathcal{V}}a$. In this case, $Cn_{S}(A) \subseteq Cn_{\mathcal{V}}(A)$, for all $A \subseteq X$.

We say that the valuation structure $(X,\mathcal{V})$ is \emph{complete} with respect to the formal system $S$ if and only if for all $A \subseteq X$ and $a \in X$ we have that: if $A \models_{\mathcal{V}}a$, then $A \vdash_{S}a$. Thus, $Cn_{\mathcal{V}}(A) \subseteq Cn_{S}(A)$, for all $A \subseteq X$. 

We say that the valuation structure $(X,\mathcal{V})$ is \emph{adequate} with respect to the formal system $S$ if and only if it is sound and complete for $S$ , i.e., for all $A \subseteq X$ and $a \in X$ we have that: $A \vdash_{S}a$ if and only if $A \models_{\mathcal{V}}a$. In this case, $Cn_{S}(A) = Cn_{\mathcal{V}}(A)$, for all $A \subseteq X$.

An important fact is that for all formal system $S = (X,\mathbb{A},\mathbb{R})$ there is a valuation structure $(X,\mathcal{V})$ which is adequate for $S$. Hence, there is no problem in considering valuations as bivaluations, i.e., functions with codomain $\{0,1\}$.\footnote{In this sense, the methodology developed here straightforwardly applies to many-valued logics assuming what is called \emph{Suszko's Thesis}: many-valued logics can be converted into bivalent logics (cf. \cite{tsuji}).}

This can be easily proved considering $THE^{*}_{S}$, the class of $S$-consistent theories, and taking as $\mathcal{V}$ the characteristic functions of elements of $THE^{*}_{S}$.

\begin{lema}
\label{lema1}
Let $S = (X,\mathbb{A},\mathbb{R})$ be a formal system and $(X,\mathcal{V})$ a valuation structure. Then, it follows:

i. If $(X,\mathcal{V})$ is complete for $S$, then $A \subseteq X$ is $S$-consistent implies that $A$ is $\mathcal{V}$-satisfiable;

ii. If $(X,\mathcal{V})$ is sound  for $S$, then $A \subseteq X$ is $\mathcal{V}$-satisfiable implies that $A$ is $S$-consistent.
\end{lema}

\begin{prova}
$i$. Consider $A \subseteq X$, $S$-consistent and assume, by \emph{reductio ad absurdum}, that $A$ is $\mathcal{V}$-unsatisfiable, i.e., $Mod_{\mathcal{V}}(A) = \emptyset$. But, then, for all $a \in X$, it holds that $A \models_{\mathcal{V}}a$. Given completeness, it follows that $A \vdash_{S}a$ for all $a \in A$, that is, $A$ is $S$-inconsistent. (contradiction!).

$ii$. Consider $A \subseteq X$, $\mathcal{V}$-satisfiable and assume, by \emph{reductio ad absurdum}, that $A$ is $S$-inconsistent, i.e, $Cn_{S}(A) = X$. Given soundness, we have that $Cn_{\mathcal{V}}(A) = X$. So, $Mod_{\mathcal{V}}(A) \subseteq Mod_{\mathcal{V}}(a)$, for all $a \in X$. Given that $A$ is $\mathcal{V}$-satisfiable, 
let $v \in Mod_{\mathcal{V}}(A)$. Then, $v \in Mod_{\mathcal{V}}(a)$ for all $a \in X$, i. e., $v(a) = 1$ for all $a \in X$, that is, $v$ is the constant function $\mathbf{1}$ (contradiction with definition of $\mathcal{V}$!).
\end{prova}

\medskip
The notion of truth which is essential to semantics appear here as an ingredient of our \emph{valuation structures} leading to the concept of semantical consequence operation $Cn_{\mathcal{V}}$. Up to now, we have settled basic usual terminology in order to separate between proof-theoretical and semantical consequence relations, and to show standard connections between them. From now on, we introduce some original concepts that constitute the main contribution of this paper.

\section{Paradeduction}

In the following, we define \emph{paradeduction} for axiomatic formal systems, though it seems
theoretically possible to realize the same work for other proof presentations such as natural deduction, sequent calculus or \emph{tableaux} systems.

Let $X$ be a set and consider a finite sequence of pairs  
$$\sigma = (A_{1},a_{1}),...,(A_{n},a_{n})$$ 
such that $A_{i} \subseteq X$ and $a_{i} \in X$, for $1 \leq i \leq n$. We introduce, based on $\sigma$, two sequences given by:
$$
\Pi_{1}\sigma := A_{1},...,A_{n} 
$$
$$
\Pi_{2}\sigma := a_1,...,a_{n}
$$
which are sequences composed by the first and second elements of the original sequence 
$\sigma$.

Consider an axiomatic formal system $S = (X,\mathbb{A},\mathbb{R})$ and a subset $A \subseteq X$. A \emph{paradeduction} in $S$ from $A$ is a finite sequence of pairs 

\begin{center}
$\sigma = (A_{1},a_{1}),...,(A_{n},a_{n})$ 
\end{center}
for $1 \leq i \leq n$, and $A_{i} \subseteq X$, $a_{i}\in X$ such that:

\begin{enumerate}
	\item $A_{i}$ is $S$-consistent, for all $1 \leq i \leq n$.
	\item For each $1 \leq i \leq n$, we have:
	
	\begin{enumerate}
		\item $a_{i} \in A$ and $A_{i} = \{a_{i}\}$ (and, therefore, $\{a_{i}\}$ is $S$-consistent) or;
		\item $a_{i} \in \mathbb{A}$ and $A_{i} = \emptyset$ or;
		\item $a_{i}$ is an immediate consequence of a set of preceding formulas $\{a_{i_{1}},...,a_{i_{k}}\}$ in the sequence $\Pi_{2}\sigma$ and $A_{i} = \bigcup_{j=1}^{k}A_{i_{j}}$ (and, also, $A_{i}$ is $S$-consistent).
	\end{enumerate}
\end{enumerate}

To illustrate, consider the following set of information:

\begin{center}
 $A=\{a \wedge b, a \rightarrow c, b \rightarrow \neg c\}$
\end{center}

If classical logic is the logic underlying $A$, it is obvious one can produce a \emph{deduction} 
such that $A \vdash c \wedge \neg c$. However, in the case of a \emph{paradeduction}, we should take
into account only consistent subsets of $A$ (in this case, $A$ is the only set which should be disregarded).
Then, a \emph{paradeduction} would allow to infer $c$ and to infer $\neg c$, 
but never their conjunction, as only consistent sets (of information) could be taken into consideration.

\begin{lema}
\label{lema2}
Let $(A_{1},a_{1}),...,(A_{n},a_{n})$ be a paradeduction in $S$ from $A \subseteq X$. Then, for $1 \leq i \leq n$, $A_{i}$ is $S$-consistent, $A_{i} \subseteq A$ and $A_{i} \vdash_{S}a_{i}$.
\end{lema}

\begin{prova}
By definition, each $A_{i}$, for $1 \leq i \leq n$, is $S$-consistent. In addition, considering the item in the definition of a paradeduction, either $A_{i} = \emptyset$, or $A_{i} = \{a_{i}\}$ or $A_{i}$ is the union of the previous $A_{k}$'s . The proof is by induction on \emph{i}. Suppose we have that $A_{i} \subseteq A$, for all $1 \leq i \leq n$.

Consider a paradeduction given by:
$$
(A_{1},a_{1}),...,(A_{n},a_{n}).
$$
For $i = 1$, we have two cases:
\begin{enumerate}
	\item $a_{1} \in A$. In this case, $A_{1} = \{a_{1}\}$ and we have, given property (I) of axiomatic formal systems, that $A_{1} \vdash_{S}a_{1}$.
	\item $a_{1} \in \mathbb{A}$. Here, $A_{1} = \emptyset$ and we have, by the definition of deduction in $S$, that $A_{1} \vdash_{S}a_{1}$.
\end{enumerate}
Suppose that $A_{i}\vdash_{S}a_{i}$, for $1 \leq i \leq k-1$, and consider $(A_{k},a_{k})$. We have, then, three cases: i. $a_{k} \in A$; ii. $a_{k} \in \mathbb{A}$ or iii. $a_{k}$ is obtained by an inference rule. In the first two cases, an argument analogous to the one above shows that $A_{k}\vdash_{S}a_{k}$. Now, for case (iii). Suppose $a_{k}$ is an immediate consequence of $\{a_{i_{1}},...,a_{i_{l}}\}$ by the application of the rule $R \in \mathbb{R}$. By induction hypothesis, we have that:
$$
\begin{array}{c}
A_{i_{1}} \vdash_{S} a_{i_{1}} \\
\vdots \\
A_{i_{l}} \vdash_{S} a_{i_{l}}
\end{array}
$$
And, moreover, $A_{k} = \bigcup_{j=1}^{l}A_{i_{j}}$ is $S$-consistent. But, then,  using properties (II) and (III) of axiomatic formal systems, the sequence $a_{i_{1}},...,a_{i_{l}},a_{k}$ is an $S$-deduction of $a_{k}$ from $A_{k} = \bigcup_{j=1}^{l}A_{i_{j}}$. That is, $A_{k}\vdash_{S}a_{k}$. 
\end{prova}

\medskip
Let $S = (X,\mathbb{A},\mathbb{R})$ be a formal system, $A \subseteq X$ and $a \in X$. We say that $a$ is $S$-\emph{paradeductible} from $A$ if there is a paradeduction in $S$, $(A_{1},a_{1}),...,(A_{n},a_{n})$, from $A$ such that $a = a_{n}$. We use, to denote this fact, the notation:
$$
A \vdash_{S}^{P}a.
$$

Now we present a result that relates the notion of paradeducibility above with a construction proposed in \cite{desouza} of a paraconsistentization functor in the category of consequence structures.

\begin{proposicao}
\label{prop1}
Let $S = (X,\mathbb{A},\mathbb{R})$ be a formal system such that $A \subseteq X$ and $a \in X$. Then, it holds that: $A \vdash_{S}^{P}a$ if and only if there is $A' \subseteq A$, $S$-consistent such that $A' \vdash_{S}a$.
\end{proposicao}

\begin{prova}
($\Rightarrow$) Suppose that $A \vdash_{S}^{P}a$. Then, there is a paradeduction in $S$ given by $(A_{1},a_{1}),...,(A_{n},a_{n})$ from $A$ such that $a = a_{n}$. By lemma \ref{lema2}, we have that $A_{i}$ is $S$-consistent, $A_{i} \subseteq A$ and $A_{i} \vdash_{S}a_{i}$, for all $1 \leq i \leq n$. Thus, there is $A_{n} \subseteq A$, $S$-consistent such that $A_{n}\vdash_{S}a$.

($\Leftarrow$) Suppose that there is  $A' \subseteq A$, $S$-consistent such that $A' \vdash_{S}a$. Consider, then, a $S$-deduction 
of $a$ from $A'$ given by:
$$
a_{1},...,a_{n} (=a).
$$
Since every subset of $A'$ is $S$-consistent (property VII), it is easy to convert the $S$-deduction above into a paradeduction in $S$:
$$
(A_{1},a_{1}),...,(A_{n},a_{n})
$$
making

\begin{equation*}
A_{i} = \left\{
\begin{array}{ll}
\{a_{i}\} & \text{If } a_{i} \in A',\\
\emptyset & \text{If } a_{i} \in \mathbb{A},\\
\bigcup_{k}A_{k} & \text{If } a_{i} \text{ is a consequence of previous formulas } \\
 & A_{k}'s \text{ associated}.
\end{array} \right.
\end{equation*}

Therefore, as $A' \subseteq A$, it follows that $A \vdash_{S}^{P}a$.                         
\end{prova}

\medskip
Consider, now, a valuation structure $(X,\mathcal{V})$. Let $A \subseteq X$ and $a \in X$. By definition, we already have that:
$$
A\models_{\mathcal{V}}a \Leftrightarrow Mod_{\mathcal{V}}(A) \subseteq Mod_{\mathcal{V}}(a).
$$

We define a relation of $\mathcal{V}$-\emph{paraconsequence} between $A$ and $a$ given by:
$$
A\models_{\mathcal{V}}^{P}a \Leftrightarrow \text{ there exists } A' \subseteq A, \;\; \mathcal{V}\text{-satisfiable such that } A'\models_{\mathcal{V}}a.
$$

We have the following result:

\begin{teorema}
Let $S = (X,\mathbb{A},\mathbb{R})$ be a formal system and $(X,\mathcal{V})$ a valuation structure adequate for $S$. Then, for all $A \subseteq X$ and $a \in X$, it is the case that:

i. $A\models_{\mathcal{V}}^{P}a$ implies that $A\vdash_{S}^{P}a$;

ii. $A\vdash_{S}^{P}a$ implies that $A\models_{\mathcal{V}}^{P}a$.
\end{teorema}

\begin{prova}
$i$. Assume that $A\models_{\mathcal{V}}^{P}a$. Then, there is $A' \subseteq A$, $\mathcal{V}$-satisfiable such that $A'\models_{\mathcal{V}}a$. By completeness, we have that $A' \vdash_{S}a$ and, by soundness and lemma \ref{lema1} (ii), $A'$ is $S$-consistent. Thus, by proposition \ref{prop1}, $A\vdash_{S}^{P}a$.

ii. Suppose that $A\vdash_{S}^{P}a$. By proposition \ref{prop1}, we have that there is $A' \subseteq A$, $S$-consistent such that $A' \vdash_{S}a$. By soundness, we have that $A' \models_{\mathcal{V}}a$ and, by completeness and lemma \ref{lema1} (i), $A'$ is $\mathcal{V}$-satisfiable. Thus, $A\models_{\mathcal{V}}^{P}a$.
\end{prova}

\medskip

The theorem above shows that if a given logic defined as a formal system and a valuation structure is sound and complete, then the paraconsistentized version of it also has these properties, because paradeduction corresponds to paraconsistentization. This is an important piece of information, especially because if one is dealing with a given theory which appears to be contradictory, then there is a guarantee that the underlying logic can be paraconsistentized in order to deal with these same contradictions as the resulting paraconsistentized logic also has desirable metalogical properties and, therefore, replacing one explosive logic by a non-explosive one is legitimate operation. We consider this topic in what follows.

\section{Paradeduction and paraconsistency}

In \cite{desouza}, a sufficient condition for the paraconsistentization functor to turn
an explosive consequence structure into a paraconsistent one is studied. Indeed, a given consequence structure should be normal, explosive, satisfying joint consistency and conjunctive properties
in order to be paraconsistentized (p.246, theorem 4.2 of  \cite{desouza}). In this section,
we show that there is an analogous argument to prove that the notion of paradeduction
paraconsistentizes axiomatic formal systems. As a natural corollary, given completeness, the same holds for the semantical notion of paraconsequence. 

Let $S = (X,\mathbb{A},\mathbb{R})$ be an axiomatic formal system. Assume that the set $X$ is endowed with an operator $\neg$ representing \emph{negation}, i.e., whenever $x \in X$, we have $\neg x \in X$. Remember that if $S = (X,\mathbb{A},\mathbb{R})$  is an axiomatic formal system and if $\vdash_{S}$ is its consequence relation, we define for all $A \subseteq X$:

\begin{center}
$Cn_{S}(A) := \{a \in X:A \vdash_{S} a \}$.
\end{center}

Thus, $Cn_{S}$ has the following properties:

i. Inclusion:  $A\subseteq Cn_{S}(A)$;

ii. Idempotency: $Cn_{S}(Cn_{S}(A)) \subseteq Cn_{S}(A)$;

iii. Monotonicity: $Cn_{S}(A) \subseteq Cn_{S}(A \cup B)$.

In this sense, using an adaptation of definition 4.1 of \cite{desouza}
for a particular consequence relation, we have for an axiomatic formal system $S$ and its associated consequence relation $\vdash_S$ that: 

(a) $\vdash_{S}$ is \emph{explosive} if and only if for all $A \subseteq X$, if there is $x \in X$ such that
$A \vdash_{S} x$ and  $A \vdash_{S} \neg x$, then $A \vdash_{S} y$, for all $y \in X$. Otherwise, $\vdash_{S}$ is
\emph{paraconsistent}. 

(b) $S$ satisfies \emph{joint consistency} if and only if there is $x \in X$ such that $\{x\}$ and $\{\neg x\}$
are both $S$-consistent while $\{x, \neg x\}$ is $S$-inconsistent. 

(c) $S$ satisfies \emph{conjunctive property} if and only if for all $x,y \in X$, there exists $z \in X$ such that $\{x, y\} \vdash_{S} u$ iff $z \vdash_{S} u$ for all $u \in X$. That is, $Cn_{S}(\{x,y\})= Cn_{S}(\{z\})$

Consider now a formal system $S = (X,\mathbb{A},\mathbb{R})$ with its consequence relation
and let $\vdash_{S}^{P}$ be the paradeduction consequence relation associated to $\vdash_{S}$.
We have the following result similar to theorem 4.2 of \cite{desouza}.

\begin{teorema}
If $S = (X,\mathbb{A},\mathbb{R})$ is an axiomatic formal system that satisfies joint consistency and
conjunctive property, and if $\vdash_{S}$ is explosive, then $\vdash_{S}^{P}$ is paraconsistent.
\end{teorema}

\begin{prova}

Since $S = (X,\mathbb{A},\mathbb{R})$ satisfies joint consistency, there is 
$a \in X$ such that $\{a\}$ and $\{\neg a\}$ are both $S$-consistent, but $Cn_{S}(\{a, \neg a\})=X$. Consider
$A=\{a, \neg a\}$. By inclusion, $A \subseteq Cn_{S}(A)$. So, by joint consistency,
$A \subseteq Cn_{S}^{P}(A)$ := $\{x \in X: A \vdash_{S}^{P} x \}$. So, $A \vdash_{S}^{P} a$ and $A \vdash_{S}^{P} \neg a$. Thus, as $S$ satisfies conjunctive property, then there exists $c \in X$ such that $Cn_{S}(\{c\})=Cn_{S}(A)=X$.
We show that $A \nvdash_{S}^{P} c$, i.e., $\vdash_{S}^{P}$ is paraconsistent. 

Notice that set $A$ has three $S$-consistent subsets: $\emptyset$, $\{a\}$ and also $\{\neg a\}$ ($\emptyset$
is $S$-consistent for $\emptyset \subseteq \{a\}$, and by monotonicity, $Cn_{S}(\emptyset) \subseteq 
Cn_{S}(\{a\}) \neq X$). We have to show that: $c \notin Cn_{S}(\emptyset)$, $c \notin Cn_{S}(\{a\})$ and $c \notin Cn_{S}(\{\neg a\})$. The argument is the same in all three cases. Let $K \in \{\emptyset, \{a\}, \{\neg a\}\}$. If $c \in Cn_{S}(K)$, then
$\{c\} \subseteq Cn_{S}(K)$. By monotonicity, $Cn_{S}(\{c\}) \subseteq Cn_{S}(Cn_{S}(K))$ and, by inclusion
and idempotency, $Cn_{S}(Cn_{S}(K))=Cn_{S}(K)$. Therefore, since $Cn_{S}(K) \neq X$, we have $X=Cn_{S}(\{c\}) \subseteq Cn_{S}(K) \neq X$ (contradiction!). So, we have shown that $A \not\vdash_S^P c$. This
implies that $\vdash_S^P$ is paraconsistent. 
\end{prova}

\section{Conclusion}

Our study offers a contribution within the domain of paraconsistent logics. Before the development of the method presented here and elsewhere (see \cite{desouza}), some forms of paraconsistentization had to be ``handcrafted", in the sense that this had to be done individually for each logic. For instance, Newton da Costa defined the consequence relation of paraclassical logic $\mathbb{P}$ as: $A \models a$ if and only if there is a $\mathbb{C}$-consistent $A'\subset A$ such that $A' \models a$.\footnote{Where $\mathbb{C}$ denotes the classical logic.} This  is a specific procedure to turn classical logic into a paraconsistent logic. Our work consists on a generalization of this idea, allowing the abstraction of the notion of classically consistent to the notion of consistent in a given logic, that is, $Cn$-consistent. It is also a generalization in the sense that da Costa's procedure was entirely semantical, while the notion of paradeduction introduced on previous sections allows a syntactical investigation of paraconsistentization, and also of the meta-linguistic properties of soundness and completeness. 

Another example would be Rescher and Manor's machinery for making inferences from inconsistent premises.\footnote{Cf. \cite{rescher-manor}.} Here, the consequence relation is not restricted to consistent subsets of the original set, but to maximal consistent subsets.\footnote{This gives rise to the notion of \emph{weak and strong consequence} A proposition $a$ is a weak consequence of a set $A$ if and only if there is at least one $A' \subseteq A$-maximal consistent, such that $a$ is a logical consequence of $A'$; and $a$ is a strong consequence of $A$ if and only if it follows from all maximal consistent subsets of $A$.} The details are not important here. What matters is that, once again, we have a particular method for dealing with inconsistent information. With respect to the work of Rescher and Manor our paper allows a second level of abstraction. Not only we can make the notion of consistency relative to a given logic, we can also replace consistency constraints for any other desired notion. Thus, Rescher and Manor's machinery becomes a particular case of our paraconsistentization method, when we change $Cn$-consistent notion for $Cn$-maximally consistent. Depending on our goal, we could even replace the notion of consistency for some desired epistemic notion, and investigate the resulting logic. So, compared to these individual ways of dealing with inconsistency that produce one paraconsistent logic at a time, our approach permits this transformation on a ``industrial scale", creating infinitely many paraconsistent logics.

In addition to that, our paper is related to belief revision theory, given that uses of consistent sets seems to be similar to some operations in belief revision theory, especially those of contraction and revision studied and developed in \cite{agm}. The main account of belief revision theory is the AGM approach. In it, theories are defined as sets closed under the consequence relation of classical logic. Thus any inconsistent theory is a trivial theory and the source of the inconsistency cannot be localized, making contraction\footnote{The contraction of a proposition $\alpha$ from a theory $\Gamma$ is a theory $\Gamma'$ that does not contain - and does not imply - $\alpha$, but it is in some sense similar to $\Gamma$.} a somewhat arbitrary operation. Furhmann\footnote{Cf. \cite{fuhrmann}.} develops a framework that is independent of classical logic. Hence, he can, in principle, accommodate inconsistent theories. But only in the sense that it allows him to localize the inconsistency and then to perform a contraction to restore the theory's consistency. Indeed, he claims that ``restoring consistency is one of the major, if not the principal reason for contracting a theory."\footnote{\cite{fuhrmann}, p. 187.} 

Paraconsistentization offers an alternative approach. When facing an inconsistent theory, we can paraconsistentize the original explosive logic, and reason with the consistent subsets of the theory. So, there is no need to contract a belief from a theory. The pragmatic advantage of one approach over the other has to be established on a case by case analysis. From a philosophical point of view, however, paraconsistentization shows that it is possible to reason with inconsistent theories without removing any of its propositions.
So, discussions regarding legal procedures as remarked in \cite{newton} or philosophical disputes can be regulated by our methodology of concentrating on consistent sets. It is also possible to explain early uses of calculus in a similar way. However, in this paper we have not explored these fruitful relationships further.

This paper showed that soundness and completeness are invariant when logics are paraconsistentized. One of the limits of our study is that it works
precisely, at first glance, only for axiomatic formal systems and valuation structures, but there are many other presentations of logics available. Thus, some problems remaining open are: To develop other forms of paraconsistentization based on the criteria presented in many of the previous works on consistent subsets, and investigate their relations. How to develop paradeduction for other proof-theoretical presentations of logics like natural deduction, sequents or \emph{tableaux}?  Other step is the investigation of problems in complexity theory. For instance, given a logic in which the satisfiability problem belongs to a given complexity class, what happens in the case of paraconsistentizing this logic? A conjecture is that complexity increases, as the set of consistent subsets has greater cardinality. However, we do not have a proof of this fact. We conjecture also that if complexity decreases, then ways of paraconsistentizing can also be used to study some problems in the domain of advanced theoretical computer science.

\subsection*{Acknowledgements}
Special thanks to two anonymous reviewers for essential comments and suggestions to improve this paper.

\end{document}